\documentclass[11pt,reqno]{amsart}
\usepackage{amsmath,amssymb}

 \makeatletter
 \evensidemargin  \oddsidemargin
 \makeatother

\begin{document}
\thispagestyle{empty}
 \setcounter{page}{1}

\begin{center}
{\Large\bf  Stability of ternary Jordan homomorphisms and
derivations associated to the generalized Jensen equation

\vskip.20in

{\small\bf  $^1$M. Eshaghi Gordji,  $^2$E. Rashidi and $^3$J. M. Rassias } \\[2mm]

{\footnotesize $^{1,2}$Department of Mathematics,
Semnan University,\\ P. O. Box 35195-363, Semnan, Iran\\
$^3$ Section of Mathematics and Informatics, Pedagogical Department
National and Capodistrian University of Athens, 4, Agamemnonos St.,
Aghia Paraskevi, Athens 15342, Greece\\
[-1mm] e-mail: {\tt madjid.eshaghi@gmail.com,
ehsanerashidi@gmail.com, jrassias@primedu.uoa.gr}}}

\end{center}
\vskip 5mm
 \noindent{\footnotesize{\bf Abstract.} In this paper, we establish the generalized
 Hyers-Ulam stability of
Jordan homomorphisms and Jordan derivations between ternary
algebras via the generalized Jensen equation
$rf(\frac{sx+ty}{r})=sf(x)+tf(y)$.

\vskip.10in
 \footnotetext {Received ????.  Revised ????. }

  \newtheorem{df}{Definition}[section]
  \newtheorem{rk}[df]{Remark}
   \newtheorem{lem}[df]{Lemma}
   \newtheorem{thm}[df]{Theorem}
   \newtheorem{pro}[df]{Proposition}
   \newtheorem{cor}[df]{Corollary}
   \newtheorem{ex}[df]{Example}

 \setcounter{section}{0}
 \numberwithin{equation}{section}

\vskip .2in

\begin{center}
\section{Introduction}
\end{center}
It seems that the stability problem of functional equation had
been first raised by Ulam [28], for what metric group $G$ is it
true that an $\varepsilon$-automorphism of $G$ is necessarily near
to a strict automorphism?.
 An answer to the above problem has been given as follows. Suppose
 $E_{1}$ and $E_{2}$ are two real Banach spaces and $f:E_{1}\longrightarrow
 E_{2}$ is a mapping. If there exists $\varepsilon\geq 0$ and $p\geq
 0,p\neq 1$ such that
 $$\|f(x+y)-f(x)-f(y)\|\leq\varepsilon(\|x\|^{p}+\|y\|^{p})\eqno(1.1)$$
for all $x,y\in E_{1}$, then there is a unique additive mapping
$h:E_{1}\longrightarrow E_{2}$ such that
$\|f(x)-h(x)\|\leq\frac{2\varepsilon \|x\|^{P}}{|2-2^{p}|}$ for
every $x\in E_{1}$. This result is called the generalized
Hyers--Ulam stability of the additive Cauchy equation
$g(x+y)=g(x)+g(y)$. Indeed, Hyers obtained the result for $p=0$.
Then Rassias [21] generalized the above result of Hyers to the case
where $0\leq p<1$. In 1994, a generalization of Th. Rassias's
theorem, was obtained by G\v avruta [7] by following the same
approach as in [21].

Following the terminology of [1], a non-empty set $G$ with a ternary
operation $[.,.,.]:G\times G\times G\longrightarrow G$ is called a
ternary groupoid and is denoted by $(G,[.,.,.])$. The ternary
groupoid $(A,[.,.,.])$ is called commutative if
$[x_{1},x_{2},x_{3}]=[x_{\sigma(1)},x_{\sigma(2)},x_{\sigma(3)}]$
for all $x_{1},x_{2},x_{3}\in A$ and all permutations $\sigma$ of
{1,2,3}.
 If a binary operation $\circ$ is defined on $G$ such that
 $[x,y,z]=x\circ y\circ z$ for all $x,y,z\in G$, we say that
 $[.,.,.]$ is derived from $\circ$. We say that $(G,[.,.,.]$ is a
 ternary semigroup if the operation $[.,.,.]$ is associative,
 i.e., if $[[x,y,z],u,v]=[x,[y,z,u],v]=[x,y,[z,u,v]$ holds for all
$x,y,z,u,v\in G$ (see[3]).
 A ternary algebra is a complex Banach space $A$, equipped with
a ternary product $(x,y,z)\longmapsto[x,y,z]$ of $A^{3}$ into $A$,
which is $\mathbb{C}$-linear in every variable, and associative in
the sense that $[x,y,[z,u,v]]=[x,[u,z,y],v]$, and satisfies
$\|[x,y,z]\|\leq\|x\|.\|y\|.\|z\|$ and $\|[x,x,x]\|=\|x\|^{3}$
(see[1]).
 If a ternary algebra $(A,[.,.,.])$ has an identity e, i.e., an
element $e\in A$ such that $x=[x,e,e]=[e,e,x]$ for all $x\in A$,
then it is routine to verify that $A$, endowed with $x\circ
y:=[x,e,y]$ and $x^{\ast}:=[e,x,e]$, is a unital algebra.
Conversely, if $(A,\circ)$ is a unital algebra, then
$[x,y,z]:=x\circ y^{\ast}\circ z$ makes $A$ into a ternary algebra.
 A $\mathbb{C}$-linear mapping $f:A\longrightarrow A$ is called:\\

$\bullet$ ternary (algebra) homomorphism if
$$f([x,y,z])=[f(x),f(y),f(z)]$$for all $x,y,z\in A$.

$\bullet$ Ternary (algebra) Jordan homomorphism if
$$f([x,x,x])=[f(x),f(x),f(x)]$$for all $x\in A$.

$\bullet$ Ternary (algebra) derivation if
$$f([x,y,z])=[f(x),y,z]+[x,f(y),z]+[x,y,f(z)]$$for all $x,y,z\in A$.

$\bullet$ Ternary (algebra) Jordan derivation if
$$f([x,x,x])=[f(x),x,x]+[x,f(x),x]+[x,x,f(x)]$$for all $x\in A$.

 A generalization of the Jensen equation is the
equation
$$rf(\frac{sx+ty}{r})=sf(x)+tf(y),$$
where $f$ is a mapping between linear spaces and $r,s,t$ are given
constant values (see [14]). A mapping $f$ satisfying $f(0)=0$ is a
solution of the Jensen functional equation if and only if satisfies
the additive Cauchy equation $f(x+y)=f(x)+f(y)$ (see [4]).
 Throughout the paper, $A$ denotes a Banach ternary algebra and
$X$ is a Banach space. Let a mapping $f$ satisfy the generalized
Jensen equation. Hence without loss of generality we can assume that
$f(0)=0$.

\vskip 5mm
\section{ternary Jordan homomorphisms}

First we obtain the stability of generalized Jensen equation by
difference Hyers sequences of [13] and [4], and then we establish
the superstability and generalized Hyers--Ulam stability of ternary
Jordan homomorphisms associated to this equation.

\begin{thm}\label{t2} $[17]$ Let $f: X\longrightarrow X$ be a mapping with
$f(0)=0$, for which there exists a function $\varphi: X\times
X\longrightarrow [0,\infty)$ satisfying
$$\phi(x,y):=\frac{1}{r}\sum(\frac {r}{s})^{-n}\varphi((\frac{r}{s})^{n}x,(\frac{r}{s})^{n}y)<\infty,$$
and
$$\|rf(\frac{sx+ty}{r})-sf(x)-tf(y)\|\leq\varphi(x,y),\eqno\hspace{1 cm}(2.1)$$
for all $ x,y \in X$. Then there exists a unique additive mapping
$h:X\longrightarrow X $given by
$h(x):=\lim_{n\longrightarrow\infty}(\frac{r}{s})^{-n}f((\frac{r}{s})^{n}x)$
satisfying $h(\frac{r}{s}x)=\frac{r}{s}h(x)$ such that
$$\|f(x)-h(x)\|\leq\phi(x,x)$$
for all $x\in X$.
\end{thm}

 Similarly we have the following theorem.
\begin{thm}\label{t2}
Let $f:X\longrightarrow X$ be a mapping with $f(0)=0$ for which
there exists a function $\varphi:X\times X\longrightarrow
[0,\infty)$ satisfying
$$\widetilde{\varphi}(x,y):=\frac{1}{s}\sum_{n=0}^{\infty}(\frac{r}{s})
^{n}\varphi((\frac{r}{s})^{-n}x,(\frac{r}{s})^{-n}y)<\infty,\eqno(2.2)$$
and
$$\|rf(\frac{sx+ty}{r})-sf(x)-tf(y)\|\leq\varphi(x,y),\eqno(2.3)$$
for all $x,y\in X$. Then there exists a unique additive mapping
$h:X\longrightarrow X$ given by $h(x):=\lim_{n\longrightarrow
\infty}(\frac{r}{s})^{n}f((\frac{r}{s})^{-n}x)$ satisfying
$h(\frac{s}{r}x)=\frac{s}{r}h(x)$ such that
$$\|f(x)-h(x)\|\leq\widetilde{\varphi}(x,x),$$
for all $x\in X$.
\end{thm}

Now we prove the superstability problem for ternary Jordan
homomorphisms as follows.
\begin{pro}\label{t2}
Let $r\neq s$ and $h:A\longrightarrow A$ be a mapping with
$h(\frac{r}{s}x)=\frac{r}{s}h(x),x\in A$ for which there exists a
function $\varphi:A^{3}\longrightarrow[0,\infty)$ satisfying
$$\lim_{n\longrightarrow
\infty}(\frac{r}{s})^{-n}\varphi((\frac{r}{s})^{n}x,(\frac{r}{s})^{n}y,(\frac{r}{s})^{n}a)=0$$
and
$$\|rh(\frac{\mu{sx}+\mu{ty}+[aaa]}{r})-\mu{s}h(x)-\mu{t}h(y)-[h(a)h(a)h(a)]\|$$
$$\leq\varphi(x,y,a),\eqno\hspace{1 cm}(2.4)$$
for all $\mu\in \mathbb C$ and all $x,y,a\in A$. Then $h$ is a
ternary Jordan homomorphisms.
\end{pro}
\begin{proof}
$h(0)=\frac{r}{s}h(0)$ and $\frac{r}{s}\neq 1$, then $h(0)=0$. Put
$\mu=1,a=0$ and replace $x,y$ by
$(\frac{r}{s})^{n}x,(\frac{r}{s})^{n}y$ in (2.4), respectively, we
get
$$\|r(\frac{r}{s})^{-n}h((\frac{r}{s})^{n}\frac{sx+ty}{r})-s(\frac{r}{s})
^{-n}h((\frac{r}{s})^{n}x)+t(\frac{r}{s})^{-n}h(\frac{r}{s})^{n}y\|$$
$$\leq(\frac{r}{s})^{-n}\varphi((\frac{r}{s})^{n}x,(\frac{r}{s})^{n}y,0).$$
If $n\longrightarrow \infty$ in above inequality, then  we conclude
that $h$ satisfies the Jensen equation. Hence $h$ is additive.
Similarly one can prove that $h(\mu{x})=\mu{h(x)}$ for all $\mu\in
\mathbb C$ all $x\in A$. Putting $x=y=0$ and replacing $a$ by
$(\frac{r}{s})^{n}a$ in (2.4), we get
$$\|rh(\frac{[aaa]}{r})-[h(a)h(a)h(a)]\|=(\frac{r}{s})^{-3n}\|rh(\frac
{[(\frac{r}{s})^{n}a((\frac{r}{s})^{n}a)(\frac{r}{s})^{n}a]}{r})$$
$$-[h((\frac{r}{s})^{n}a)h((\frac{r}{s})^{n}a)h((\frac{r}{s})^{n})]\|$$
$$\leq(\frac{r}{s})^{-3n}\varphi(0,0,(\frac{r}{s})^{n}a)$$
$$\leq(\frac{r}{s})^{-n}\varphi(0,0,(\frac{r}{s})^{n}a),$$
for all $a\in A$. The right hand side of above inequality tends to
zero as $n\longrightarrow \infty$, so that
$$h([aaa])=rh(\frac{[aaa]}{r})=[h(a)h(a)h(a)]$$
for all $a\in A$. Hence,  $h$ is a ternary Jordan homomorphism.
\end{proof}
\begin{thm}\label{t2}
Let $f:A\longrightarrow A$ be a mapping such that $f(0)=0$, for
which there exists a function $\varphi:A^{3}\longrightarrow
[0,\infty)$ such that
$$\widetilde{\varphi}(x,y,a):=\frac{1}{r}\sum_{j=0}^{\infty}(\frac{r}{s})
^{-j}\varphi((\frac{r}{s})^{j}x,(\frac{r}{s})^{j}y,(\frac{r}{s})^{j}a)\leq\infty,\eqno\hspace{1
cm}(2.5)$$ and
$$\|rf(\frac{\mu{sx}+\mu{ty}+[aaa]}{r})-\mu{sf(x)}-\mu{tf(y)}-[f(a)f(a)f(a)]\|$$
$$\leq\varphi(x,y,a),\eqno\hspace{1 cm}(2.6)$$
for all $\mu\in \mathbb{T}^{1}=\{Z\in\mathbb{C}:|z|=1\}$ and all
$a\in A$. Then there exists a unique ternary Jordan homomorphism
$h:A\longrightarrow B$ such that
$$\|f(x)-h(x)\|\leq\widetilde{\varphi}(x,x,0)\eqno\hspace{1 cm}(2.7)$$
for all $x\in A$
\end{thm}
\begin{proof}
Put $a=0$ and $\mu=1$ in (2.6). By using theorem 2.1, there is a
unique additive mapping $h:A\longrightarrow B$ defined by
$h(x)=\lim_{n\longrightarrow
\infty}(\frac{r}{s})^{-n}f((\frac{r}{s})^{n}x),$ and satisfying
(2.7) for all $x\in X$.
 Let $\mu\in\mathbb{T}^{1}$. Replace $x$ by $(\frac{r}{s})^{n+1}x$
 and put $y=0$ in (2.6), then
 $$\|f((\frac{r}{s})^{n}\mu{x})-\mu(\frac{r}{s})^{-1}f((\frac{r}{s})^{n+1}x)\|
 \leq\frac{1}{r}\varphi((\frac{r}{s})^{n+1}x,0,0),$$
 for all $x\in A$. It follows from the inequality
  $\|rf(x)-sf(\frac{r}{s}x)\|\leq\varphi(\frac{r}{s}x,0,0)$ that
$$\|f((\frac{r}{s})^{n}x)-(\frac{r}{s})^{-1}f((\frac{r}{s})^{n+1}x)\|\leq\frac{1}{s}\varphi((\frac{r}{s})^{n+1}x,0,0).$$
Hence,
$$(\frac{r}{s})^{-n}f((\frac{r}{s})^{n}\mu{x})-\mu(\frac{r}{s})^{-n}f((\frac{r}{s})^{n}x)\|$$
$$\leq(\frac{r}{s})^{-n}\|f((\frac{r}{s})^{n}\mu{x})-\mu(\frac{r}{s})^{-1}f((\frac{r}{s})^{n+1}x)\|$$
$$+(\frac{r}{s})^{-n}\|\mu{f((\frac{r}{s})^{n}x)}-\mu(\frac{r}{s})^{-1}f((\frac{r}{s})^{n+1}x)\|$$
$$\leq2(\frac{r}{s})^{-n}(\frac{1}{r}\varphi((\frac{r}{s})^{n+1}x,0,0)$$
for all $x\in A$. The right hand side of above inequality tends to
zero as $n\longrightarrow\infty$. By using (2.3), we get
$$h(\mu{x})=\lim_{n\longrightarrow\infty}(\frac{r}{s})^{-n}f((\frac{r}
{s})^{n}\mu{x})=\lim_{n\longrightarrow\infty}(\mu(\frac{r}{s})^{-n}f((\frac{r}{s})^{n}x))=\mu{h(x)},$$
for all $x\in A$.  Let $\lambda\in\mathbb{C} (\lambda\neq 0)$ and
let $M$ be a natural number greater than $4|\lambda|$. Then
$|\frac{\lambda}{M}|<\frac{1}{4}<1-\frac{2}{3}=\frac{1}{3}$. By
theorem 1 of [10], there exist three numbers
$\mu_{1},\mu_{2},\mu_{3}\in\mathbb{T}^{1}$ such that
$3\frac{\lambda}{M}=\mu_{1}+\mu_{2}+\mu_{3}$. By the additivity of
$h$ we get $h(\frac{1}{3})=\frac{1}{3}h(x)$ for all $x\in A$.
Therefore
$$h(\lambda{x})=h(\frac{M}{3}\cdot3\cdot\frac{\lambda}{M}x)=Mh(\frac{1}
{3}\cdot3\cdot\frac{\lambda}{M}x)=\frac{M}{3}h(3\cdot\frac{\lambda}{M}x)$$
$$=\frac{M}{3}h(\mu_{1}x+\mu_{2}x+\mu_{3}x)=\frac{M}{3}(h(\mu_{1}x)+h(\mu_{2}x)+h(\mu{3}x))$$
$$\frac{M}{3}(\mu_{1}+\mu_{2}+\mu_{3})h(x)=\frac{M}{3}\cdot(3\cdot\frac{\lambda}{M})h(x)=\lambda h(x)$$
for all $x\in A$. So $h$ is $\mathbb{C}$-linear. Set $\mu=1$ and
$x=y=0$, and replace $a$ by $(\frac{r}{s})^{n}a$, respectively, in
(2.6). Then
$$(\frac{r}{s})^{-3n}\|rf((\frac{r}{s})^{3n}\frac{[aaa]}{r})-[f((\frac{r}{s})
^{n}a)f((\frac{r}{s})^{n}a)f((\frac{r}{s})^{n}a)]\|$$
$$\leq(\frac{r}{s})^{-3n}\varphi(0,0,(\frac{r}{s})^{n}a),$$ for all
$a\in A$. Then by applying the continuity of the ternary product
$(x,y,z)\longmapsto[xyz]$ we deduce
$$h([aaa])=rh(\frac{1}{r}[aaa])$$
$$=\lim_{n\longrightarrow\infty}(\frac{r}{s})^{-3n}rf((\frac{r}{s})^{3n}\frac{[aaa]}{r})$$
$$=\lim_{n\longrightarrow\infty}[(\frac{r}{s})^{-n}f((\frac{r}{s})^{n}a)(\frac{r}{s})^{-n}
f((\frac{r}{s})^{n}a)(\frac{r}{s})^{-n}f((\frac{r}{s})^{n}a)]$$
$$=[h(a)h(a)h(a)],$$
for all $a\in A$. Thus $h$ is a ternary Jordan homomorphism
satisfying the required inequality.
\end{proof}
\begin{cor}\label{t2}
Let $f:A\longrightarrow A$ be a mapping such that $f(0)=0$, for
which there exists $\varepsilon\geq 0$ and $p<1$ such that
$$\|rf(\frac{\mu{sx}+\mu{ty}+[aaa]}{r})-\mu{sf(x)}-\mu{tf(y)}-[f(a)f(a)f(a)]\|$$
                                           $$\leq\varepsilon(\|x\|^{p}+\|y\|^{p}+\|a\|^{p}),$$
for all $\mu\in\mathbb{T}^{1}$, and all $x,y,a\in A$. Then there
exists a unique ternary Jordan homomorphism $h:A\longrightarrow A$
such that
$$\|f(x)-h(x)\|\leq\frac{2{r}^{-p}\varepsilon{\|x\|}^{p}}{{r}^{1-p}-{s}^{1-p}},$$
for all $x\in A$.
\end{cor}
\begin{proof}
Define $\varphi(x,y,a)=\varepsilon(\|x\|^{p}+\|y\|^{p}+\|a\|^{p})$,
and apply theorem 2.4.
\end{proof}
{\bf Remark 2.6.} When $p>1$,one may use the same techniques used
in the proof of theorem 2.4 to get a result similar to corollary
2.5.
\begin{cor}\label{t2}
Let $A$ be linearly spanned by a set $S\subseteq A$ and
$f:A\longrightarrow A$ be a mapping with $f(0)=0$ satisfying
$f((\frac{r}{s})^{2n}[s_{1}s_{2}a])=[f((\frac{r}{s})^{n}s_{1})f((\frac{r}{s})^{n}s_{2})f(a)]$
for all sufficiently large positive integers $n$, and all
$s_{1},s_{2}\in S,a\in A$. Suppose that
$$\|rh({\mu{sx}+\mu{ty}+[aaa]}{r})-\mu{sh(x)}-\mu{th(y)}-[h(a)h(a)h(a)]\|$$
                                                           $$\leq\varphi(x,y,a),$$
for all $\mu\in\mathbb{T}^{1}$ and all $x,y\in A$. Then there
exists a unique ternary Jordan homomorphism $h:A\longrightarrow A$
satisfying (2.7) for all $x\in A$.
\end{cor}
\begin{proof}
By the same argument as in the proof of theorem 2.4, there exists
a unique linear mapping $h:A\longrightarrow A$ given by
$$h(x):=\lim_{n\longrightarrow\infty}(\frac{r}{s})^{-n}f((\frac{r}{s})^{n}x),\hspace{1
cm}(x\in A),$$ such that
$$\|f(x)-h(x)\|\leq\widetilde{\varphi}(x,x,0),$$
for all $x\in A$. Then
$$h([s_{1}s_{2}a])=\lim_{n
\longrightarrow\infty}(\frac{r}{s})^{-2n}f([((\frac{r}{s})^{n}s_{1})((\frac{r}{s})^{n}s_{2})a])$$
$$=\lim_{n\longrightarrow\infty}[(\frac{r}{s})^{-n}f((\frac{r}{s})^{n}s_{1})(\frac{r}
{s})^{-n}f((\frac{r}{s})^{n}s_{2})f(a)]$$
$$=[h(s_{1})h(s_{2})h(s_{3})].$$
By linearity of $h$, we have
$$h([aaa])=[h(a)h(a)h(a)]$$
for all $a\in A$. Therefore
$(\frac{r}{s})^{n}h([aaa])=h([aa(\frac{r}{s})^{n}a)])=[h(a)h(a)f((\frac{r}{s})^{n}a)]$,
and so
$$h([aaa])=\lim_{n\longrightarrow\infty}(\frac{r}{s})^{-n}[h(a)h(a)f((\frac{r}{s})^{n}a)]$$
$$=[h(a)h(a)\lim_{n\longrightarrow\infty}(\frac{r}{s})^{-n}f((\frac{r}{s})^{n}a)]$$
$$=[h(a)h(a)h(a)],$$
for all $a\in A$.
\end{proof}
\begin{thm}\label{t2}
Suppose that $f:A\longrightarrow A$ is mapping with $f(0)=0$ for
which there exists a function $\varphi:A^{3}\longrightarrow
[0,\infty)$ fulfilling $(2.5)$, and $(2.6)$ for $\mu=1,\textbf{i}$
and all $x\in A$. Then there exists a unique ternary Jordan
homomorphism $h:A\longrightarrow A$ such that
$$\|f(x)-h(x)\|\leq\widetilde{\varphi}(x,x,0),$$
for all $x\in A$.
\end{thm}
\begin{proof}
Put $a=0$ and $\mu=1$ in (2.6). By the same argument as theorem
2.4 we infer that there exists a unique additive mapping
$h:A\longrightarrow A$ given by
$$h(x):=\lim_{n\longrightarrow\infty}(\frac{r}{s})^{-n}f((\frac{r}{s})^{n}x),$$
and satisfying (2.7) for all $x\in A$. By the same reasoning as in
the proof of the main theorem of [21], the mapping $h$ is
$\mathbb{R}$-linear.
 Replace $x$ by $(\frac{r}{s})^{n}x$, $y$ by $0$,  respectively, and put
 $\mu=\textbf{i}$ and $a=0$ in (2.6). Then
 $$(\frac{r}{s})^{-n}\|\frac{r}{s}f(\textbf{i}(\frac{r}{s})^{n-1}x)-\textbf{i}f
 ((\frac{r}{s})^{n}x)\|\leq\frac{1}{s}(\frac{r}{s})^{-n}\varphi((\frac{r}{s})^{n}x,0,0),$$
 for all $x\in A$. The right hand side of above inequality tends to zero as
 $n\longrightarrow\infty$, hence
 $$h(\textbf{i}x)=\lim_{n\longrightarrow\infty}(\frac{r}{s})^{-n+1}f((\frac{r}{s})
 ^{n-1}\textbf{i}x)=\lim_{n\longrightarrow\infty}\textbf{i}(\frac{r}{s})^{-n}f((\frac{r}{s})^{n}x)=\textbf{i}h(x),$$
 for all $x\in A$.
 For every $\lambda\in\mathbb{C}$ we can write
 $\lambda=\alpha_{1}+\textbf{i}\alpha_{2}$ in which
 $\alpha_{1},\alpha_{2}\in\mathbb{R}$. Therefore
 \begin{align*}h(\lambda{x})&=h(\alpha_{1}x+\textbf{i}\alpha_{2}x)=\alpha_{1}h(x)+\alpha_{2}h(\textbf{i}x)\\
 &=\alpha_{1}h(x)+\textbf{i}\alpha_{2}h(x)=(\alpha_{1}+\textbf{i}\alpha_{2})h(x)\\
 &=\lambda{h(x)},
 \end{align*}
for all $x\in A$. Hence,  $h$ is $\mathbb{C}$-linear.
\end{proof}

\vskip 5mm
\section{ternary Jordan derivations}
In this section, we establish the superstability and generalized
Hyers--Ulam stability of Jordan derivations  associated to the
generalized Jensen equation.
\begin{pro}\label{t2}
Let $r\neq s$ and $h:A\longrightarrow A$ be a mapping with
$h(\frac{r}{s}x)=\frac{r}{s}h(x),$ for all  $x\in A$,  for which
there exists a function $\varphi:A^{3}\longrightarrow[0,\infty)$
satisfying
$$\lim_{n\longrightarrow
\infty}(\frac{r}{s})^{-n}\varphi((\frac{r}{s})^{n}x,(\frac{r}{s})^{n}y,(\frac{r}{s})^{n}a)=0$$
and
$$\|rh(\frac{\mu{sx}+\mu{ty}+[aaa]}{r})-\mu{s}h(x)-\mu{t}h(y)-[h(a),a,a]-[a,h(a),a]-[a,a,h(a)]\|$$
$$\leq\varphi(x,y,a),\eqno(3.1)$$
for all $\mu\in \mathbb C$ and all $x,y,a\in A$. Then $h$ is a
ternary Jordan derivation.
\end{pro}
\begin{proof} It follows from
$h(0)=\frac{r}{s}h(0)$ and $\frac{r}{s}\neq 1$ that  $h(0)=0$. Put
$\mu=1,a=0$ and replace $x,y$ by
$(\frac{r}{s})^{n}x,(\frac{r}{s})^{n}y$, respectively,  in (3.1), we
get
$$\|r(\frac{r}{s})^{-n}h((\frac{r}{s})^{n}\frac{sx+ty}{r})-s(\frac{r}{s})^{-n}
h((\frac{r}{s})^{n}x)+t(\frac{r}{s})^{-n}h(\frac{r}{s})^{n}y\|$$
$$\leq(\frac{r}{s})^{-n}\varphi((\frac{r}{s})^{n}x,(\frac{r}{s})^{n}y,0).$$
If $n\longrightarrow \infty$, in above inequality, then $h$
satisfies the Jensen equation. Hence $h$ is additive. Similarly one
can prove that $h(\mu{x})=\mu{h(x)}$ for all $\mu\in \mathbb C$ all
$x\in A$. Putting $x=y=0$ and replacing $a$ by $(\frac{r}{s})^{n}a$
in (3.1), we get

\begin{align*}
\|rh(\frac{[a,a,a]}{r})-[h(a),a,a]&-[a,h(a),a]-[a,a,h(a)]\|=(\frac{r}{s})^{-3n}\|rh
(\frac{[(\frac{r}{s})^{n}a,(\frac{r}{s})^{n}a,(\frac{r}{s})^{n}a]}{r})\\
&-[h((\frac{r}{s})^{n}a),(\frac{r}{s})^{n}a,(\frac{r}{s})^{n}a]-[(\frac{r}{s})^{n}a,
h((\frac{r}{s})^{n}a),(\frac{r}{s})^{n}a]\\
&-[(\frac{r}{s})^{n}a,(\frac{r}{s})^{n}a,h((\frac{r}{s})^{n}a)\|\leq(\frac{r}{s})^{-3n}\varphi(0,0,(\frac{r}{s})^{n}a)\\
&\leq(\frac{r}{s})^{-n}\varphi(0,0,(\frac{r}{s})^{n}a),
\end{align*}
for all $a\in A$. The right hand side of above inequality tends to
zero as $n\longrightarrow \infty$, so that
$$h([a,a,a])=rh(\frac{[a,a,a]}{r})=[h(a),a,a]+[a,h(a),a]+[a,a,h(a)]$$
thus $h$ is a ternary Jordan derivation.
\end{proof}
\begin{thm}\label{t2}
Let $f:A\longrightarrow A$ be a mapping such that $f(0)=0$, for
which there exists a function $\varphi:A^{3}\longrightarrow
[0,\infty)$ such that
$$\widetilde{\varphi}(x,y,a):=\frac{1}{r}\sum_{j=0}^{\infty}(\frac{r}{s})^{-j}\varphi((\frac{
r}{s})^{j}x,(\frac{r}{s})^{j}y,(\frac{r}{s})^{j}a)\leq\infty,\eqno(3.2)$$
and
$$\|rf(\frac{\mu{sx}+\mu{ty}+[aaa]}{r})-\mu{s}f(x)-\mu{t}f(y)-[f(a),a,a]-[a,f(a),a]-[a,a,f(a)]\|$$
$$\leq\varphi(x,y,a),\eqno(3.3)$$
for all $\mu\in \mathbb{T}^{1}=\{Z\in\mathbb{C}:|z|=1\}$ and all
$a\in A$. Then there exists a unique ternary Jordan derivation
$h:A\longrightarrow A$ such that
$$\|f(x)-h(x)\|\leq\widetilde{\varphi}(x,x,0)\eqno(3.4)$$
for all $x\in A$
\end{thm}
\begin{proof}
Put $a=0$ and $\mu=1$ in (3.3). By using theorem 2.1, there is a
unique additive mapping $h:A\longrightarrow B$ defined by
$h(x)=\lim_{n\longrightarrow
\infty}(\frac{r}{s})^{-n}f((\frac{r}{s})^{n}x),$ and satisfying
(3.4) for all $x\in X$.
 Let $\mu\in\mathbb{T}^{1}$. Replace $x$ by $(\frac{r}{s})^{n+1}x$
 and put $y=0$ in (2.6), then
 $$\|f((\frac{r}{s})^{n}\mu{x})-\mu(\frac{r}{s})^{-1}f((\frac{r}{s})^{n+1}x)\|\leq\frac{1}{r}
 \varphi((\frac{r}{s})^{n+1}x,0,0),$$
 for all $x\in A$. It follows from
  $\|rf(x)-sf(\frac{r}{s}x)\|\leq\varphi(\frac{r}{s}x,0,0)$ that
$$\|f((\frac{r}{s})^{n}x)-(\frac{r}{s})^{-1}f((\frac{r}{s})^{n+1}x)\|\leq\frac{1}{s}\varphi((\frac{r}{s})^{n+1}x,0,0).$$
Hence,
\begin{align*}(\frac{r}{s})^{-n}&f((\frac{r}{s})^{n}\mu{x})-\mu(\frac{r}{s})^{-n}f((\frac{r}{s})^{n}x)\|\\
&\leq(\frac{r}{s})^{-n}\|f((\frac{r}{s})^{n}\mu{x})-\mu(\frac{r}{s})^{-1}f((\frac{r}{s})^{n+1}x)\|\\
&+(\frac{r}{s})^{-n}\|\mu{f((\frac{r}{s})^{n}x)}-\mu(\frac{r}{s})^{-1}f((\frac{r}{s})^{n+1}x)\|\\
&\leq2(\frac{r}{s})^{-n}(\frac{1}{r}\varphi((\frac{r}{s})^{n+1}x,0,0)
\end{align*}
for all $x\in A$. The right hand side of above inequality tends to
zero as $n\longrightarrow\infty$. By using (2.3), we get
$$h(\mu{x})=\lim_{n\longrightarrow\infty}(\frac{r}{s})^{-n}f((\frac{r}{s})^{n}\mu{x})
=\lim_{n\longrightarrow\infty}(\mu(\frac{r}{s})^{-n}f((\frac{r}{s})^{n}x))=\mu{h(x)},$$
for all $x\in A$.  We have, $h(0x)=0=0h(x)$. Let
$\lambda\in\mathbb{C} (\lambda\neq 0)$ and let $M$ be a positive
integer greater than $4|\lambda|$. Then
$|\frac{\lambda}{M}|<\frac{1}{4}<1-\frac{2}{3}=\frac{1}{3}$. By
theorem 1 of [10], there exist three numbers
$\mu_{1},\mu_{2},\mu_{3}\in\mathbb{T}^{1}$ such that
$3\frac{\lambda}{M}=\mu_{1}+\mu_{2}+\mu_{3}$. By the additivity of
$h$ we get $h(\frac{1}{3})=\frac{1}{3}h(x)$ for all $x\in A$.
Therefore
\begin{align*}h(\lambda{x})&=h(\frac{M}{3}\cdot3\cdot\frac{\lambda}{M}x)=Mh(\frac{1}{3}\cdot3\cdot
\frac{\lambda}{M}x)=\frac{M}{3}h(3\cdot\frac{\lambda}{M}x)\\
&=\frac{M}{3}h(\mu_{1}x+\mu_{2}x+\mu_{3}x)=\frac{M}{3}(h(\mu_{1}x)+h(\mu_{2}x)+h(\mu_{3}x))\\
&\frac{M}{3}(\mu_{1}+\mu_{2}+\mu_{3})h(x)=\frac{M}{3}\cdot3\cdot\frac{\lambda}{M}h(x)=\lambda
h(x) \end{align*} for all $x\in A$. So that $h$ is
$\mathbb{C}$-linear. Set $\mu=1$ and $x=y=0$, and replace $a$ by
$(\frac{r}{s})^{n}a$, respectively, in (3.3). Then
$$(\frac{r}{s})^{-3n}\|rf((\frac{r}{s})^{3n}\frac{[a,a,a]}{r})-[f((\frac{r}{s})^{n}a),(\frac{r}{s})^{n}a,
(\frac{r}{s})^{n}a]$$
$$-[(\frac{r}{s})^{n}a,f((\frac{r}{s})^{n}a),(\frac{r}{s})^{n}a]-[(\frac{r}{s})^{n}a,(\frac{r}{s})^{n}a,
f((\frac{r}{s})^{n}a)]\|$$
$$\leq(\frac{r}{s})^{-3n}\varphi(0,0,(\frac{r}{s})^{n}a),$$ for all
$a\in A$. Then by applying the continuity of the ternary product
$(x,y,z)\longmapsto[xyz]$ we deduce that
\begin{align*}h([a,a,a])&=rh(\frac{1}{r}[a,a,a])=\lim_{n\longrightarrow\infty}
(\frac{r}{s})^{-3n}rf((\frac{r}{s})^{3n}\frac{[a,a,a]}{r})\\
&=\lim_{n\longrightarrow\infty}(\frac{r}{s})^{-3n}([f((\frac{r}{s})^{n}a),(\frac{r}{s})^{n}a,
(\frac{r}{s})^{n}a]\\
&+[(\frac{r}{s})^{n}a,f((\frac{r}{s})^{n}a),(\frac{r}{s})^{n}a]+[(\frac{r}{s})
^{n}a,(\frac{r}{s})^{n}a,f((\frac{r}{s})^{n}a)])\\
&=\lim_{n\longrightarrow\infty}([(\frac{r}{s})^{-n}f((\frac{r}{s})^{n}a),a,a]+[a,(\frac{r}{s})
^{-n}f((\frac{r}{s})^{n}a),a]\\
&+[a,a,(\frac{r}{s})^{-n}f((\frac{r}{s})^{n}a)])=[h(a),a,a]+[a,h(a),a]+[a,a,h(a)]
\end{align*}
for all $a\in A$. Hence, $h$ is a ternary Jordan derivation
satisfying the required inequality.
\end{proof}
\begin{cor}\label{t2}
Let $f:A\longrightarrow A$ be a mapping such that $f(0)=0$, for
which there exist $\varepsilon\geq 0$ and $p<1$ such that
$$\|rf(\frac{\mu{sx}+\mu{ty}+[a,a,a]}{r})-\mu{sf(x)}-\mu{tf(y)}-[f(a),a,a]+[a,f(a),a]+[a,a,f(a)]\|$$
                                           $$\leq\varepsilon(\|x\|^{p}+\|y\|^{p}+\|a\|^{3p}),$$
for all $\mu\in\mathbb{T}^{1}$, and all $x,y,a\in A$. Then there
exists a unique ternary Jordan derivation $h:A\longrightarrow A$
such that
$$\|f(x)-h(x)\|\leq\frac{2{r}^{-p}\varepsilon{\|x\|}^{p}}{{r}^{1-p}-{s}^{1-p}},$$
for all $x\in A$.
\end{cor}
\begin{proof}
Define
$\varphi(x,y,a)=\varepsilon(\|x\|^{p}+\|y\|^{p}+\|a\|^{3p})$, and
apply theorem $3.2$.
\end{proof}
\begin{cor}\label{t2}
Let $A$ be linearly spanned by a set $S\subseteq A$ and
$f:A\longrightarrow A$ be a mapping with $f(0)=0$ satisfying
$$f((\frac{r}{s})^{2n}[s_{1},s_{2},a])=[f((\frac{r}{s})^{n}s_{1}),(\frac{r}{s})
^{n}s_{2},a]+[(\frac{r}{s})^{n}s_{1},f((\frac{r}{s})^{n}s_{2}),a]+[(\frac{r}{s})^{n}s_{1},(\frac{r}{s})^{n}s_{2},f(a)]$$
for all sufficiently large positive integers $n$, and all
$s_{1},s_{2}\in S,a\in A$.suppose that
$$\|rh(\frac{\mu{sx}+\mu{ty}+[a,a,a]}{r})-\mu{sh(x)}-\mu{th(y)}-[h(a),a,a]-[a,h(a),a]-[a,a,h(a)]\|$$
                                                           $$\leq\varphi(x,y,a),$$
for all $\mu\in\mathbb{T}^{1}$ and all $x,y\in A$. Then there
exists a unique ternary Jordan derivation $h:A\longrightarrow A$
satisfying (3.4) for all $x\in A$.
\end{cor}
\begin{proof}
By the same argument as in the proof of theorem $3.2$, there
exists a unique linear mapping $h:A\longrightarrow A$ given by
$$h(x):=\lim_{n\longrightarrow\infty}(\frac{r}{s})^{-n}f((\frac{r}{s})^{n}x),\hspace{1
cm}(x\in A).$$ Such that
$$\|f(x)-h(x)\|\leq\widetilde{\varphi}(x,x,0),$$
for all $x\in A$. Then \begin{align*}h([s_{1},s_{2},a])&=\lim_{n
\longrightarrow\infty}(\frac{r}{s})^{-2n}f([((\frac{r}{s})^{n}s_{1})((\frac{r}{s})^{n}s_{2})a])\\
&=\lim_{n\longrightarrow\infty}([(\frac{r}{s})^{-n}f((\frac{r}{s})^{n}s_{1}),(\frac{r}{s})^{-n}
(\frac{r}{s})^{n}s_{2},a]\\
&+[(\frac{r}{s})^{-n}(\frac{r}{s})^{n}s_{1},(\frac{r}{s})^{-n}f((\frac{r}{s})^{n}s_{2}),a]\\
&+[(\frac{r}{s})^{-n}(\frac{r}{s})^{n}s_{1},(\frac{r}{s})^{-n}(\frac{r}{s})^{n}s_{2},f(a)]\\
&=[h(s_{1}),s_{2},a]+[s_{1},h(s_{2}),a]+[s_{1},s_{2},h(a)].\end{align*}
By the linearity of $h$ we have
$$h([a,a,a])=[h(a),a,a]+[a,h(a),a]+[a,a,h(a)]$$
for all $a\in A$. Therefore
$$(\frac{r}{s})^{n}h([a,a,a])=h([a,a,(\frac{r}{s})^{n}a)])=[h(a),a,(\frac{r}{s})^{n}a]+
[a,h(a),(\frac{r}{s})^{n}a]+[a,a,f((\frac{r}{s})^{n}a)],$$ and so
\begin{align*}h([aaa])&=\lim_{n\longrightarrow\infty}(\frac{r}{s})^{-n}([h(a),a,(\frac{r}{s})^{n}a]+
[a,h(a),(\frac{r}{s})^{n}a]+[a,a,f((\frac{r}{s})^{n}a)])\\
&=[h(a),a,(\frac{r}{s})^{-n}(\frac{r}{s})^{n}a]+[a,h(a),(\frac{r}{s})^{-n}(\frac{r}
{s})^{n}a]+[a,a,\lim_{n\longrightarrow\infty}(\frac{r}{s})^{-n}f((\frac{r}{s})^{n}a)]\\
&=[h(a),a,a]+[a,h(a),a]+[a,a,h(a)], \end{align*} for all $a\in A$.
\end{proof}
{\small

}
\end{document}